\documentclass[11pt]{amsart}

\usepackage{amssymb,amsthm,amsmath}
\usepackage[numbers,sort&compress]{natbib}
\usepackage[margin=1in]{geometry}
\usepackage{color}
\usepackage{graphicx}
\usepackage{tikz}
\usepackage[colorlinks,
linkcolor=red,
anchorcolor=blue,
citecolor=blue
]{hyperref}

\def\a{\alpha}

\def\e{\varepsilon}

\let\newpf\proof \let\proof\relax 
\newenvironment{pf}{\newpf[\proofname]}{\qed\endtrivlist}

\newcommand{\ba}{\overline{A}}

\def\be{\begin{equation}}
\def\ee{\end{equation}}

\def\ba{{\begin{align}}}
\def\ea{{\end{align}}}

\def\bm{\begin{matrix}}
\def\em{\end{matrix}}

\def\a{{\alpha}}

\def\0{{\mathbf 0}}

\newtheorem{Theorem}{Theorem}[section]
\newtheorem{Lemma}{Lemma}[section]
\newtheorem{Proposition}{Proposition}[section]
\newtheorem{Corollary}{Corollary}[section]
\newtheorem{Remark}{Remark}[section]

\newtheorem{Definition}{Definition}[section]

\numberwithin{equation}{section}

\theoremstyle{definition}

\newtheorem{definition}{Definition}[section]

\renewcommand{\mod}{\operatorname{mod}}

\newcommand{\C}{{\mathbb C}}

\newcommand{\N}{{\mathbb N}}
\newcommand{\Q}{{\mathbb Q}}
\newcommand{\R}{{\mathbb R}}
\newcommand{\T}{{\mathbb T}}

\newcommand{\Z}{{\mathbb Z}}

\def\B0{{\bold{0}}}

\catcode`\@=12

\def\Empty{}
\newcommand\oplabel[1]{
  \def\OpArg{#1} \ifx \OpArg\Empty {} \else
    \label{#1}
  \fi}

\newcommand{\comm}[1]{}
\newcommand{\comment}[1]{}

\begin{document}

\title[]{On the almost reducibility conjecture}

\author{Lingrui Ge}
\address{Beijing International Center for Mathematical Research, Peking University, Beijing, China}
\email{gelingrui@bicmr.pku.edu.cn}

\begin{abstract}
Avila's Almost Reducibility Conjecture (ARC) is a powerful statement linking purely analytic and
dynamical properties of analytic one frequency $SL(2,\R)$ cocycles. It is also a fundamental tool in the study of spectral theory
of analytic one-frequency Schr\"odinger operators, with many striking
consequencies, allowing to give a detailed characterization of the
subcritical region. Here we give a proof, completely different from
Avila's, for the important case of Schr\"odinger cocycles with trigonometric polynomial potentials and
non-exponentially approximated frequencies, allowing,
in particular, to obtain all the desired spectral consequences in this case.
\end{abstract}

\maketitle
\section{Introduction}
Consider the following analytic one-frequency Schr\"odinger operator,
\begin{equation}
(H_{V,\alpha,\theta}u)_n=u_{n-1}+u_{n+1}+V(\theta+n\alpha)u_n,\ \ n\in\Z,
\end{equation}
where $\alpha\in\R\backslash\Q$ is the frequency, $\theta\in\R$ is the phase and $V\in C_h^\omega(\T,\R)$ \footnote{Let $F$ be a bounded analytic (possibly matrix valued) function defined on $ \{ \theta |  | \Im \theta |< h \}$,
$\| F\| _h= \sup_{ | \Im \theta |< h } \| F(\theta)\| $. $C^\omega_{h}(\T,*)$ denotes the
set of all these $*$-valued functions ($*$ will usually denote $\R$, $sl(2,\R)$ or
$SL(2,\R)$). Denote by $C^{\omega}(\T,\R)$  the union $\cup_{h>0}C_h^{\omega}(\T,\R)$.} is the potential.

%Avila's \cite{Aglobal} established the elegant global theory and he considered the complexified Lyapunov exponent,
For $E\in\R,$ the complexified Lyapunov exponent is defined by
\begin{equation}\label{le}
L_\e(E)=\lim\limits_{n\rightarrow\infty}\frac{1}{n}\int_{\T} A_E(\theta+i\e+(n-1)\alpha)\cdots A_E(\theta+i\e+\alpha)A_E(\theta+i\e)d\theta,
\end{equation}
where
$$
A_E(\theta+i\e)=\begin{pmatrix}E-V(\theta+i\e)&-1\\ 1&0\end{pmatrix}
$$
and $(\alpha,A_E)$ is called the Schr\"odinger cocycle.

Similarly, for general analytic $SL(2,\R)$ cocycles, using \eqref{le} with $A_E$ replaced by $A$ also defines the complexified Lyapunov exponent $L_\e(A).$

Avila's global theory of one-frequency Schr\"odinger operators
\cite{Aglobal} is based on the classification of corresponding
Schr\"odinger cocycles into three regimes, with the terminology
inspired by the almost Mathieu operator (see e.g. \cite{AvilaJito,Jito1999})
\begin{enumerate}
\item Subcritical regime: $L_0(A)=0$ and $\omega(A)=0$;
\item Critical regime: $L_0(A)=0$ and $\omega(A)>0$;
\item Supercritical regime: $L_0(A)>0$ and $\omega(A)>0$,
\end{enumerate}
where the acceleration is defined by
$$
\omega(A):=\lim\limits_{\e\rightarrow 0^+}\frac{L_\e(A)-L_0(A)}{2\pi \e}.
$$

Avila showed in \cite{Aglobal} that critical regime is atypical. It is
expected that the other two (open) regimes lead to many universal features
of corresponding energy regimes.

For the subcritical regime, key to this universality is Avila's almost
reducibility conjecture,  saying that for all analytic one-frequency
$SL(2,\R)$ cocycles
\\

\noindent {\bf Conjecture:} $(\alpha,A)$ is almost reducible if and
only if $\omega(A)=0$.\\

The ARC serves as a fundamental tool to solve various spectral
problems, for example, the non-critical Dry Ten Martini Problem
\cite{ayz}, the Aubry-Andre-Jitomirskaya phase transition conjecture
\cite{ayz1,gyzh}.

The almost reducibility conjecture was first formulated in
\cite{AvilaJito}, before even the acceleration was introduced in
\cite{Aglobal}. It was first proved for  operators with
non-perturbatively small analytic potentials (in the
Bourgain-Jitomirskaya regime \cite{bj}) by Avila-Jitomirskaya
\cite{AvilaJito} and Hou-You \cite{HouYou} provided $\beta(\alpha)=0$,
where
$$
\beta(\alpha)=\limsup\limits_{k\rightarrow\infty}-\frac{\ln \|k\alpha\|_{\R/\Z}}{|k|}\geq 0
$$
is the exponential rate of approximation of
$\alpha$ by the rationals and $\|x\|_{\R/\Z}=\inf_{j\in\Z}|x-j|$.

If $\beta(\alpha)>0$, the
complete ARC was proved by Avila \cite{A3}. The complete proof of ARC
for $\beta(\alpha)=0$ was given by Avila around 2011, who gave a
series of lectures with a detailed proof at the Fields Institute,
Toronto, in 2011, as well as at an Oberwolfach Arbeitsgemeinschaft on
Quasiperiodic Schr\"odinger operators in 2012. Avila's approach, is purely dynamical and completely different.
 While Avila's proof,  to appear in \cite{A2}, has been
fully verified by the community, it has been  desirable to have a
published article, providing a detailed proof also for the (full
measure) case $\beta(\alpha)=0$ , to fully  justify, in
particular, various important spectral applications.

In this paper, we give a different proof of the ARC for the important
case of
Schr\"odinger cocycles with $\beta(\alpha)=0$ and trigonometric polynomial potentials  by establishing almost localization of the dual
operator which is of independent interest.  While we emphasize that Avila's dynamical proof holds for
the entire class of $SL(2,\R)$ cocycles, our result is
sufficient for the spectral  applications in the trigonometric polynomial case.

One can divide the proof of the ARC for Schr\"odinger cocycles into two cases depending on vanishing/nonvanishing of the Lyapunov exponent.
\begin{enumerate}
\item $L_\e(E)=L_0(E)>0$ for all $|\e|<\e_0$;
\item $L_\e(E)=L_0(E)=0$ for all $|\e|<\e_0$.
\end{enumerate}
ARC for Case (1) follows directly from the uniform hyperbolicity. ARC for Case (2) is also known as the subcritical conjecture.

\begin{Definition}
We say $(\alpha,A_E)$ is subcritical if for some $\e_0>0$,
$$
L_\e(E)=0, \ \ \forall |\e|<\e_0.
$$
In particular, the subcritical radius is defined by
$$
h(E)=\sup_{\e>0}\{\e: L_\e(E)=0\}.
$$
\end{Definition}
\begin{Theorem}\label{main 2}
For any trigonometric polynomial potential $V(x)$, assume $(\alpha,
A_E)$ is subcritical and $\beta(\alpha)=0$, then $(\alpha, A_E)$ is almost
reducible in the strongest sense. More precisely, for any $0<r<h(E)$, there exists $B_n\in C_{r}^\omega(\T,PSL(2,\R))$ and $R\in SL(2,\R)$ such that
$$
\|B_n(\theta+\alpha)^{-1}A_E(\theta)B_n(\theta)-R\|_{r}\rightarrow 0.
$$
\end{Theorem}
\begin{Remark}
We obtain almost reducibility in a uniform band $\{z:|\Im z|<r\}$ for any  $0<r<h(E)$ which is crucial for sharp spectral applications.
\end{Remark}

\section{Preliminaries}

\subsection{Continued Fraction Expansion} Let $\alpha\in (0,1)\backslash\Q$, $a_0:=0$ and $\alpha_0:=\alpha$. Inductively, for $k\geq 1$, we define
$$
a_k:=[\alpha_{k-1}^{-1}], \  \ \alpha_k=\alpha_{k-1}^{-1}-a_k.
$$
Let $p_0:=0$, $p_1:=1$, $q_0:=1$, $q_1:=a_1$. We define inductively $p_k:=a_kp_{k-1}+p_{k-2}$, $q_k:=a_kq_{k-1}+q_{k-2}$. Then $q_n$ are the  denominators of the best rational approximations of $\alpha$ since we have $\|k\alpha\|_{\R/\Z}\geq \|q_{n-1}\alpha\|_{\R/\Z}$ for all $k$ satisfying $\forall 1\leq k< q_n$, and
$$
\frac{1}{2q_{n+1}}\leq \|q_n\alpha\|_{\R/\Z}\leq \frac{1}{q_{n+1}}.
$$

\subsection{Cocycles and Lyapunov exponents}
Let ${\rm Sp}(2d,\C)$ be the set of all $2d\times 2d$ complex symplectic matrices and $\alpha$ be irrational. Given $A \in C^0(\T,{\rm Sp}(2d,\C))$, we define the complex one-frequency symplectic cocycle $(\alpha,A)$ by:
$$
(\alpha,A)\colon \left\{
\begin{array}{rcl}
\T \times \C^{2d} &\to& \T \times \C^{2d}\\[1mm]
(\theta,v) &\mapsto& (\theta+\alpha,A(\theta)\cdot v)
\end{array}
\right.  .
$$
The iterates of $(\alpha,A)$ are of the form $(\alpha,A)^n=(n\alpha,A_n)$, where
$$
A_n(\theta):=
\left\{\begin{array}{l l}
A(\theta+(n-1)\alpha) \cdots A(\theta+\alpha) A(\theta),  & n\geq 0\\[1mm]
A^{-1}(\theta+n\alpha) A^{-1}(\theta+(n+1)\alpha) \cdots A^{-1}(\theta), & n <0
\end{array}\right.    .
$$
We denote by $L_1(A)\geq L_2(A)\geq...\geq L_d(A)$ the non-negative Lyapunov exponents of $(\alpha,A)$ repeatedly according to their multiplicities, i.e.,
$$
L_k(A)=\lim\limits_{n\rightarrow\infty}\frac{1}{n}\int_{\T}\ln\sigma_k(A_n(\theta))d\theta,
$$
where  $\sigma_1(A_n)\geq...\geq \sigma_d(A_n)$ denote its singular values (eigenvalues of $\sqrt{A_n^*A_n}$) with modulus greater than or equal to one. Since the k-th exterior product $\Lambda^kA_n$ satisfies $\sigma_1(\Lambda^kA_n)=\|\Lambda^kA_n\|$, $L^k(A)=\sum_{j=1}^kL_j(A)$ satisfies

$$
L^k(A)=\lim\limits_{n\rightarrow \infty}\frac{1}{n}\int_{\T}\ln\|\Lambda^kA_n(\theta)\|d\theta.
$$

A typical example is the quasiperiodic finite-range cocycle $(\alpha,L_{E,V}^{2\cos})$ where $V(\theta)=\sum\limits_{k=-d}^dV_ke^{2\pi ik\theta}$ with $V_k=\overline{V_{-k}}$ and
\begin{align*}\tiny
L_{E,V}^{2\cos}(\theta)=\frac{1}{V_d}
\begin{pmatrix}
-V_{d-1}&\cdots&-V_1&E-V_0-2\cos2\pi(\theta)&-V_{-1}&\cdots&-V_{-d+1}&-V_{-d}\\
V_d& \\
& &  \\
& & & \\
\\
\\
& & &\ddots&\\
\\
\\
& & & & \\
& & & & & \\
& & & & & &V_{d}&
\end{pmatrix}
\end{align*}
Let
$$
C=\begin{pmatrix}
V_d&\cdots&V_1\\
0&\ddots&\vdots\\
0&0&V_d
\end{pmatrix},\ \ \Omega=\begin{pmatrix}0&-C^*\\
C&0\end{pmatrix}.
$$
One can check that $L^{2\cos}_{E,V}(\theta)$ is complex symplectic with respect to $\Omega$ for any $E\in\R$. In this case, we denote the non-negative Lyapunov exponents by $\gamma_i(E)=L_i(L_{E,V}^{2\cos})$ for $1\leq i\leq d$ for short.

\subsection{The rotation number}
We consider the one-frequency $SL(2,\R)$-cocycle $(\alpha,A)$ where $A \in C^0(\T, {\rm SL}(2, \R))$. Assume that $A$ is homotopic to the identity, then $(\alpha, A)$ induces the projective skew-product $F_A\colon \T \times \mathbb{S}^1 \to \T \times \mathbb{S}^1$
$$
F_A(x,w):=\left(x+\a,\, \frac{A(x) \cdot w}{|A(x) \cdot w|}\right),
$$
which is also homotopic to the identity. Lift $F_A$ to a map $\widetilde{F}_A\colon \T \times \R \to \T \times \R$ of the form $\widetilde{F}_A(x,y)=(x+\alpha,y+\psi_x(y))$, where for every $x \in \T$, $\psi_x$ is $\Z$-periodic.
%Let us denote by $\pi\colon \T^d \times \R\to \T^d \times \mathbb{S}^1$ the projection $(x,y) \mapsto (x,e^{2 \pi\mathrm{i} y})$. Then,
%$$
%F_A \circ \pi = \pi \circ \widetilde{F}_A.
%$$
The map $\psi\colon\T \times \R  \to \R$ is called a {\it lift} of $A$. Let $\mu$ be any probability measure on $\T \times \R$ which is invariant by $\widetilde{F}_A$, and whose projection on the first coordinate is given by Lebesgue measure.
The number
$$
\rho(\alpha,A):=\int_{\T^d \times \R} \psi_x(y)\ d\mu(x,y) \ {\rm mod} \ \Z
$$
 depends  neither on the lift $\psi$ nor on the measure $\mu$, and is called the \textit{fibered rotation number} of $(\alpha,A)$ (see \cite{H,JM} for more details).

Given $\theta\in\T$, let $
R_\theta:=
\begin{pmatrix}
\cos2 \pi\theta & -\sin2\pi\theta\\
\sin2\pi\theta & \cos2\pi\theta
\end{pmatrix}$.
If $A\colon \T\to{\rm PSL}(2,\R)$ is homotopic to $x \mapsto R_{nx/2}$ for some $n\in\Z$,
then we call $n$ the {\it degree} of $A$ and denote it by $\deg A$.
The fibered rotation number is invariant under real conjugacies which are homotopic to the identity. More generally, if $(\alpha,A_1)$ is conjugated to $(\alpha, A_2)$, i.e., $B(x+\alpha)^{-1}A_1(x)B(x)=A_2(x)$, for some $B \colon \T\to{\rm PSL}(2,\R)$ with $\deg{B}=n$, then
\begin{equation}\label{rotation number}
\rho(\alpha, A_1)= \rho(\alpha, A_2)+ \frac{\langle n,\alpha \rangle}2.
\end{equation}

In particular, for quasiperiodic Schr\"odinger cocycle $(\alpha,A_E)$ where
\begin{equation}\label{S}
A_E(x)=\begin{pmatrix}
E-V(x)&-1\\ 1&0
\end{pmatrix},
\end{equation}
we denote by $\rho(E)=\rho(\alpha,A_E)$.

\subsection{Aubry duality}

Consider the fiber direct integral,
$$
\mathcal{H}:=\int_{\T}^{\bigoplus}\ell^2(\Z)dx,
$$
which, as usual, is defined as the space of $\ell^2(\Z)$-valued, $L^2$-functions over the measure space $(\T,dx)$.  The extensions of the
Sch\"odinger operators  and their finite-range duals to  $\mathcal{H}$ are given in terms of their direct integrals, which we now define.
Let $\alpha\in\T$ be fixed. Interpreting $H_{V,\alpha,x}$ as fibers of the decomposable operator,
$$
H_{V,\alpha}:=\int_{\T}^{\bigoplus}H_{V,\alpha,x}dx,
$$
then the family $\{H_{V,\alpha,x}\}_{x\in\T}$ naturally induces an operator on the space $\mathcal{H}$, i.e. ,
$$
(H_{V,\alpha} \Psi)(x,n)= \Psi(x,n+1)+ \Psi(x,n-1) +  V(x+n\alpha) \Psi(x,n),
$$
where $V$ is a 1-periodic real trigonometric polynomial.

Similarly,  the direct integral of finite-range operator  $L_{V,\alpha,\theta}$,
denote as $L_{V,\alpha}$, is given by
$$
(L_{V,\alpha}  \Psi)(\theta,n)=  \sum\limits_{k=-d}^{d} V_k \Psi(\theta,n+k)+  2\cos2\pi (\theta+n\alpha) \Psi(\theta,n),
$$
where $V_k$ is the $k$-th Fourier coefficient of $V(x)$.

Let  $U$ be the following operator on $\mathcal{H}:$
$$
(\mathcal{U}\phi)(\eta,m)=\sum_{n\in\Z}\int_{\T}e^{2\pi imx}e^{2\pi in(m\alpha+\eta)}\phi(x,n)dx.
$$
Then direct computations show that $U$ is unitary and satisfies
$$U H_{V,\alpha} U^{-1}=L_{V,\alpha}.$$
Thus the quasiperiodic  finite-range operator  $L_{V,\alpha,\theta}$ is called the dual operator of $H_{V,\alpha,x}$ \cite{gjls}.

\subsection{Quantitative global theory}
We let $V(\theta)=\sum_{k=-d}^dV_ke^{2\pi ik\theta}$ with $\overline{V}_k=V_{-k}$. We formulate the dual characterization of the subcritical Schr\"odinger cocycle by quantitative global theory developed in \cite{gjyz}.  Assume $(\alpha,A_E)$ is subcritical, $h(E)$ is the subcritical radius and  $\gamma_d(E)$ is the minimal Lyapunov exponent of the dual finite-range cocycle $(\alpha,L_{E,V}^{2\cos})$, then

\begin{Theorem}[\cite{gjyz}]\label{sub}
For any $\alpha\in\R\backslash\Q$ and $E\in \R$, we have
$$
h(E)=\frac{\gamma_d(E)}{2\pi}>0.
$$
\end{Theorem}
Thus the fact that $(\alpha, A_E)$ is subcritical implies that $\gamma_d(E)>0$, in other words, $(\alpha,L_{E,V}^{2\cos})$ is supercritical.
\section{Upper bound of the numerator of the green's function}
In \cite{AvilaJito}, Avila-Jitomirskaya developed a method to prove  ARC for quasiperiodic operators with non-perturbatively small analytic potentials. Their proof is based on further development of the localization method introduced by Bourgain-Jitomirskaya \cite{bj} and they derive almost reducibility from almost localization of the dual long-range operators.  In their argument, a key assumption is that the potential of the dual long-range operator is large enough so that they have an estimate of the upper/lower bound of the numerator of the Green's function. The crucial part of our proof of ARC is to develop a new method to remove the largeness assumption and give sharp upper bound of the numerator of the Green's function based on dual Lyapunov exponents.

We are concerned  with the following quasiperiodic finite-range operators,
\begin{equation}\label{fo}
(L_{V,\alpha,\theta}u)_n=\sum\limits_{k=-d}^d V_k u_{n+k}+2\cos2\pi(\theta+n\alpha),\ \ n\in\Z.
\end{equation}

\subsection{The transfer matrix and determinant of truncated operators} Given an interval $I\subset \Z$, let $R_I$ be the projection of $\Z$ to $I$.

It is well known that in the Schr\"odinger case, if we denote
$$
P_k(V,\theta,E)=\det{R_{[0,k-1]}(H_{V,\alpha,\theta}-E)R^*_{[0,k-1]}}
$$
the determinant of the truncated operator with Dirichlet boundary at $-1$ and $k$, and we denote the $k$-th transfer matrix by $(A_E)_k(\theta)=A_E(\theta+(k-1)\alpha)\cdots A_E(\theta+\alpha)A_E(\theta)$, then we have the following beautiful relation between the truncated operator and  elements of the transfer matrix,
$$
(A_E)_k(\theta)=\begin{pmatrix}P_k(V,\theta,E)&-P_{k-1}(V,\theta+\alpha,E)\\
P_{k-1}(V,\theta,E)&-P_{k-2}(V,\theta+\alpha,E)\end{pmatrix}.
$$

In the following, we obtain similar results for the finite-range operators \eqref{fo}.  For any $1\leq m\leq d$,  let $\Lambda^m(L^{2\cos}_{E,V})_k(\theta)$ be the $m$-th wedge of $(L^{2\cos}_{E,V})_k(\theta)$, where
$$
(L^{2\cos}_{E,V})_k(\theta)=L^{2\cos}_{E,V}(\theta+(k-1)\alpha)\cdots L^{2\cos}_{E,V}(\theta+\alpha)L^{2\cos}_{E,V}(\theta).
$$

Let $\delta_1$, $\delta_2$, $\cdots$, $\delta_{2d}$ be the canonical orthonormal basis of $\R^{2d}$ and $-d\leq i_1,i_2,\cdots,i_m,j_1,j_2,\cdots,j_m\leq d-1$, we denote
\begin{align*}
(Q_k)_{i_1,\cdots,i_m}^{j_1,\cdots,j_m}(V,\theta,E)=\langle \delta_{i_1}\wedge\cdots\wedge\delta_{i_m},\Lambda^m(L^{2\cos}_{E,V})_k(\theta)\delta_{j_1}\wedge\cdots\wedge\delta_{j_m}\rangle.
\end{align*}
\begin{Proposition}\label{p1}
$(Q_k)_{i_1,\cdots,i_m}^{j_1,\cdots,j_m}(V,\theta,E)$ is zero  for some $E\in\C$ if and only if $L_{V,\alpha,\theta}u=Eu$ has a solution obeying the following boundary conditions
$$
u(\ell)=0, \ \ \ell\neq j_1,\cdots,j_m, \ \ -d\leq \ell\leq d-1,
$$
$$
u(k+\ell)=0, \ \ \ell= i_1,\cdots,i_m, \ \ -d\leq \ell\leq d-1.
$$
\end{Proposition}
\begin{pf}
Note that $L_{V,\alpha,\theta}u=Eu$ if and only if
\begin{equation}\label{ss1}
\begin{pmatrix}
u(k+d-1)\\ \vdots\\ u(k-d)
\end{pmatrix}=(L_{E,V}^{2\cos})_k(\theta)\begin{pmatrix}
u(d-1)\\ \vdots\\ u(-d)
\end{pmatrix}
\end{equation}
Proposition \ref{p1} follows from the definition of wedge and \eqref{ss1}.
\end{pf}
\begin{Proposition}\label{p2}
$L_{V,\alpha,\theta}u=Eu$ has a solution obeying the following boundary condition
$$
u(\ell)=0, \ \ \ell\neq j_1,\cdots,j_m, \ \ -d\leq \ell\leq d-1,
$$
$$
u(k+\ell)=0, \ \ \ell= i_1,\cdots,i_m, \ \ -d\leq \ell\leq d-1,
$$
if and only if
$$
\det{R_{[0,k-1]}(L_{V,\alpha,\theta}-E)R^*_{[d,k-d-1]\cup\{j_1,\cdots,j_m\}\cup[k-d,k+d-1]\backslash\{k+i_1,\cdots,k+i_m\}}}=0.
$$
\end{Proposition}
\begin{pf}
Simple linear algebra conclusion.
\end{pf}
Propostion \ref{p1} and \ref{p2} imply the following theorem.
\begin{Theorem}\label{th1}
We have
\begin{align*}
(Q_k)_{i_1,\cdots,i_m}^{j_1,\cdots,j_m}(V,\theta,E)=C(d)V_d^{-k}\det{R_{[0,k-1]}(L_{V,\alpha,\theta}-E)R^*_{[d,k-d-1]\cup\{j_1,\cdots,j_m\}\cup[k-d,k+d-1]\backslash\{k+i_1,\cdots,k+i_m\}}}
\end{align*}
for some $C(d)$ not depending on $k$.
\end{Theorem}
\begin{pf}
Note that for any fixed $V$ and $\theta$,
$(Q_k)_{i_1,\cdots,i_m}^{j_1,\cdots,j_m}(V,\theta,E)$ and
$$
\det{R_{[0,k-1]}(L_{V,\alpha,\theta}-E)R^*_{[d,k-d-1]\cup\{j_1,\cdots,j_m\}\cup[k-d,k+d-1]\backslash\{k+i_1,\cdots,k+i_m\}}}
$$
are both polynomials of $E$. If all the zeros of them are simple, by Proposition \ref{p1} and \ref{p2}, one has
\begin{align*}
(Q_k)_{i_1,\cdots,i_m}^{j_1,\cdots,j_m}(V,\theta,E)=C(d)V_d^{-k}\det{R^*_{[0,k-1]}(L_{V,\alpha,\theta}-E)R_{[d,k-1]\cup\{j_1,\cdots,j_m\}\cup[k,k+d-1]\backslash\{k+i_1,\cdots,k+i_m\}}}.
\end{align*}
Otherwise, one can assume
$$
(Q_k)_{i_1,\cdots,i_m}^{j_1,\cdots,j_m}(V,\theta,E)=\sum\limits_{j=0}^{p}a_j(V,\theta)E^j
$$
where $a_j(V,\theta)$ are polynomials of $V_{-d},\cdots,V_d$. Note that $Q_k:=(Q_k)_{i_1,\cdots,i_m}^{j_1,\cdots,j_m}(V,\theta,E)$ has multiple roots if and only if
$$
{\rm Res}(Q_k,Q_k')=0,
$$
where ${\rm Res}(Q_k,Q_k')$ is a polynomial of $(V_{-d},\cdots,V_{d})$, thus if we denote by
$$
\mathcal{Z}=\{(V_{-d},\cdots,V_{d}):{\rm Res}(Q_k,Q_k')(V_{-d},\cdots,V_{d})=0\},
$$
then $|\mathcal{Z}|=0$ where $|\cdot|$ is the Lebesgue measure in $\R^{2d+1}$. Thus there exists a sequence
$$
(V_{-d}^j,\cdots,V_{d}^j)\rightarrow (V_{-d},\cdots,V_d),\ \ j\rightarrow \infty
$$
such that the zeroes of $(Q_k)_{i_1,\cdots,i_m}^{j_1,\cdots,j_m}(V^j,\theta,E)$ and
$$
\det{R_{[0,k-1]}(L_{V^j,\alpha,\theta}-E)R^*_{[d,k-d-1]\cup\{j_1,\cdots,j_m\}\cup[k-d,k+d-1]\backslash\{k+i_1,\cdots,k+i_m\}}}
$$
are simple where $V^j(x)=\sum_{k=-d}^d V_k^j e^{2\pi ikx}$. Finally, both $(Q_k)_{i_1,\cdots,i_m}^{j_1,\cdots,j_m}(V,\theta,E)$ and
$$
\det{R_{[0,k-1]}(L_{V,\alpha,\theta}-E)R^*_{[d,k-d-1]\cup\{j_1,\cdots,j_m\}\cup[k-d,k+d-1]\backslash\{k+i_1,\cdots,k+i_m\}}}
$$
are continuous in $V$, we thus finish the proof.
\end{pf}
\begin{Corollary}\label{mainco}
For any $\e>0$, there exists $C(\e,d)$, such that
\begin{align*}
&\sup_{\theta\in\T}\left|\det{R_{[0,k-1]}(L_{V,\alpha,\theta}-E)R^*_{[d,k-d-1]\cup\{j_1,\cdots,j_m\}\cup[k-d,k+d-1]\backslash\{k+i_1,\cdots,k+i_m\}}}\right|\\
\leq &C(\e,d)e^{(\gamma_1(E)+\cdots+\gamma_m(E)+\ln |V_d|+\e)|k|}.
\end{align*}
\end{Corollary}
\begin{pf}
It follows from Theorem \ref{th1} and the upper semi-continuity of the Lyapunov exponents.
\end{pf}
\subsection{The representation of Green's function}
Let $I=[x_1,x_2]$ and $x_2-x_1+1=kd$.  For any $x,y\in I$, we denote
$$
G_I(x,y)=\langle\delta_x, (R_IL_{V,\alpha,\theta}R^*_I-E)^{-1}\delta_y\rangle,
$$
then for any $u$ such that $L_{V,\alpha,\theta}u=Eu$, we have
$$
u(x)=-\sum\limits_{y=x_1}^{x_1+d-1}\sum\limits_{k=y+1-x_1}^dG_I(x,y)V_{k}u(y-k)-\sum\limits_{y=x_2-d+1}^{x_2}\sum\limits_{k=-d}^{y-x_2-1}G_I(x,y)V_{k}u(y-k).
$$
Moreover, by Cramer's rule we have
$$
G_I(x,y)=\frac{\mu_{x,y}}{P_{kd}(V,\theta+x_1\alpha,E)}
$$
where $\mu_{x,y}$ is the corresponding minor and $P_{kd}(V,\theta,E)=\det{R_{[0,kd-1]}(L_{V,\alpha,\theta}-E)R^*_{[0,kd-1-1]}}$.

In the following, we give a new representation formula of $\mu_{x,y}$. Let $x_1\leq x\leq x_1+d-1$. Note that by the definition of $G_I(x,y)$, we have
$$
\mu_{x,y}=(-1)^{x+y}\det{R_{[x_1,x_2]\backslash\{y\}}(L_{V,\alpha,\theta}-E)R^*_{[x_1,x_2]\backslash\{x\}}}.
$$
As a corollary of Lemma \ref{lat}, we have
\begin{Corollary}\label{representation}
Assume $k\geq k_0+10>20$, $x\in [x_1,x_1+d-1]$ and $y\in [x_1+k_0d,x_1+(k_0+1)d-1]$, we have
$$
|\mu_{x,y}|\leq \sum\limits_{\sigma}\sum\limits_{\tau}|\det\left({R_{[x_1+k_0d,x_1+(k_0+2)d-1]\backslash\{y\}}(L_{V,\alpha,\theta}-E)R^*_{\sigma\cup\tau}}\right)\mu_{\sigma}(x,y)\mu^{\tau}(x,y)|
$$
where the sum is taken over all ordered subset $\sigma=(\sigma_1,\cdots,\sigma_{d})$ of $[x_1+(k_0+1)d,x_1+(k_0+3)d-1]$ and all ordered subset $\tau=(\tau_1,\cdots,\tau_{d-1})$ of $[x_1+(k_0-1)d+1,x_1+(k_0+1)d-1]$. Moreover,
$$
\mu_{\sigma}(x,y)=\det{R_{[x_1+(k_0+2)d,x_2]}(L_{V,\alpha,\theta}-E)R^*_{[x_1+(k_0+1)d,x_2]\backslash\sigma}},
$$
$$
\mu^{\tau}(x,y)=\det{R_{[x_1,x_1+k_0d-1]}(L_{V,\alpha,\theta}-E)R^*_{[x_1,x_1+(k_0+1)d-1]\backslash\{\tau,x\}}}.
$$
\end{Corollary}
Corollary \ref{representation} and Corollary \ref{mainco} imply the following upper bound of $\mu_{x,y}$.
\begin{Lemma}
Assume $\e>0$, $x\in [x_1,x_1+d-1]$ and $y\in [x_1+k_0d,x_1+(k_0+1)d-1]$, if $k\geq k_0+10>20$, we have
$$
|\mu_{x,y}|\leq C(\e,d)e^{(\sum\limits_{i=1}^{d-1}\gamma_i+\ln |V_d|+\e)|y-x_1|+(\sum\limits_{i=1}^{d}\gamma_i+\ln |V_d|+\e)|y-x_2|}.
$$
\end{Lemma}
\begin{pf}
For  all ordered subset $\sigma=(\sigma_1,\cdots,\sigma_{d})$ of $[x_1+(k_0+1)d,x_1+(k_0+3)d-1]$, by  Corollary \ref{mainco}, we have
\begin{align*}
|\mu_\sigma(x,y)|&=|\det{R_{[x_1+(k_0+2)d,x_2]}(L_{V,\alpha,\theta}-E)R^*_{[x_1+(k_0+1)d,x_2]\backslash\sigma}}|\\
&\leq C(\e,d)e^{(\sum\limits_{i=1}^{d}\gamma_i(E)+|\ln V_d|+\e)(k-k_0-2)d}\leq C(\e,d)e^{(\sum\limits_{i=1}^{d}\gamma_i(E)+|\ln V_d|+\e)|y-x_2|}.
\end{align*}
For all ordered subset $\tau=(\sigma_1,\cdots,\sigma_{d-1})$ of $[x_1+(k_0-1)d+1,x_1+(k_0+1)d-1]$, by  Corollary \ref{mainco}, we have
\begin{align*}
|\mu^{\tau}(x,y)|&=|\det{R_{[x_1,x_1+k_0d-1]}(L_{V,\alpha,\theta}-E)R^*_{[x_1,x_1+(k_0+1)d-1]\backslash\{\tau,x\}}}|\\
&\leq C(\e,d) e^{(\sum\limits_{i=1}^{d-1}\gamma_i(E)+|\ln V_d|+\e)k_0d}\leq C(\e,d)e^{(\sum\limits_{i=1}^{d-1}\gamma_i(E)+|\ln V_d|+\e)|y-x_1|}.
\end{align*}
Thus
$$
|\mu_{x,y}|\leq C(\e,d)e^{(\sum\limits_{i=1}^{d-1}\gamma_i(E)+\ln |V_d|+\e)|y-x_1|+(\sum\limits_{i=1}^{d}\gamma_i(E)+\ln |V_d|+\e)|y-x_2|}.
$$
\end{pf}

\section{A linear algebra trick}
In this section, we involve a simple linear algebra trick to calculate the minor of a block diagonal matrix. We consider the following block diagonal matrix,
$$
M=\begin{pmatrix}
B&A\\
A^*&B&A\\
&\ddots&\ddots&\ddots\\
&&A^*&B&A\\
&&&A^*&B
\end{pmatrix}
$$
where $A,B$ are both $d\times d$ matrices and the size of $M$ is $kd\times kd$.  Let $M(i,j)$ be the $ij$-th minor of $M$ and $P_I:[1,kd]\rightarrow I$ be the projection where $I\subset[1,kd]$.
\begin{Lemma}\label{lat}
Assume $k\geq 10+k_0\geq 20$, $k_0d+1\leq i\leq (k_0+1)d$ and $(k-1)d+1\leq j\leq kd$, we have
$$
|M(i,j)|\leq \sum\limits_{\sigma}\sum\limits_{\tau}|\det\left({P_{[(k_0-1)d+1,(k_0+1)d]\backslash\{i\}}MP^*_{\sigma\cup\tau}}\right)\mu_{\sigma}(i,j)\mu^{\tau}(i,j)|
$$
where the sum is taken over all ordered subset $\sigma=(\sigma_1,\cdots,\sigma_{d})$ of $[(k_0-2)d+1,k_0d]$ and all ordered subset $\tau=(\tau_1,\cdots,\tau_{d-1})$ of $[k_0d+1,(k_0+2)d]$. Moreover,
$$
\mu_{\sigma}(i,j)=\det{P_{[1,(k_0-1)d]}MP^*_{[1,k_0d]\backslash\sigma}},
$$
$$
\mu^{\tau}(i,j)=\det{P_{[(k_0+1)d+1,kd]}MP^*_{[k_0d+1,kd]\backslash\{\tau,j\}}}.
$$
\end{Lemma}
\begin{pf}
Note that
$$
M(i,j)=(-1)^{i+j}\det{P_{[1,kd]\backslash\{i\}}MP^*_{[1,kd]\backslash\{j\}}}.
$$
Expanding $P_{[1,kd]\backslash\{i\}}MP^*_{[1,kd]\backslash\{j\}}$ by the $\{(k_0-1)d+1,\cdots, (k_0+1)d-1\}$-rows, then we have
\begin{align}\label{we1}
|M(i,j)|\leq \sum\limits_{\gamma}\left|\det\left({P_{[(k_0-1)d+1,(k_0+1)d]\backslash\{i\}}MP^*_{\gamma}}\right)\det\left({P_{[1,kd]\backslash[(k_0-1)d+1,(k_0+1)d]}MP^*_{[1,kd]\backslash\{\gamma\cup\{j\}\}}}\right)\right|
\end{align}
where the sum is taken over all ordered subset $\gamma=(\gamma_1,\cdots,\gamma_{2d-1})$ of $[1,kd]$.
\begin{Proposition}\label{ppp1}
If $\gamma_{2d-1}\geq (k_0+2)d+1$ or $\gamma_{1}\leq (k_0-2)d$, then
$$
\det\left({P_{[(k_0-1)d+1,(k_0+1)d]\backslash\{i\}}MP^*_{\gamma}}\right)=0,
$$
where $\gamma=(\gamma_1,\cdots,\gamma_{2d-1})$.
\end{Proposition}
\begin{pf}
Note that
$$
P_{[(k_0-1)d+1,(k_0+1)d]}M=\begin{pmatrix}&&A^*&B&A&&\\
&&&A^*&B&A&\end{pmatrix},
$$
thus if $\gamma_{2d-1}\geq (k_0+2)d+1$ or $\gamma_{1}\leq (k_0-2)d$, then
$$
P_{[(k_0-1)d+1,(k_0+1)d]\backslash\{i\}}MP^*_{\gamma}
$$
contains at least one zero column.
\end{pf}
It follows Proposition \ref{ppp1} that   the sum in \eqref{we1} is taken over all ordered subset $\gamma=(\gamma_1,\cdots,\gamma_{2d-1})$ of $[(k_0-2)d+1,(k_0+2)d]$.
\begin{Proposition}\label{ppp2}
If $\#\left\{\{\gamma_j\}_{j=1}^{2d-1}\cap[(k_0-2)d+1,k_0d]\right\}>d$ or $\#\left\{\{\gamma_j\}_{j=1}^{2d-1}\cap[k_0d+1,(k_0+2)d]\right\}>d-1$, then
$$
\det\left({P_{[1,kd]\backslash[(k_0-1)d+1,(k_0+1)d]}MP^*_{[1,kd]\backslash\{\gamma\cup\{j\}\}}}\right)=0.
$$
\end{Proposition}
\begin{pf}
We distinguish two cases:

Case I: $\#\left\{\{\gamma_j\}_{j=1}^{2d-1}\cap[(k_0-2)d+1,k_0d]\right\}>d$, we have
$$
P_{[1,kd]\backslash[(k_0-1)d+1,(k_0+1)d]}MP^*_{[1,kd]\backslash\{\gamma\cup\{j\}\}}=\begin{pmatrix}A_1&0\\ A_2&A_3
\end{pmatrix}
$$
where the last column of $A_1$ is zero. Hence
$$
\det\left({P_{[1,kd]\backslash[(k_0-1)d+1,(k_0+1)d]}MP^*_{[1,kd]\backslash\{\gamma\cup\{j\}\}}}\right)=\det{A_1}\det{A_3}=0.
$$

Case II: $\#\left\{\{\gamma_j\}_{j=1}^{2d-1}\cap[k_0d+1,(k_0+2)d]\right\}>d-1$, we have
$$
P_{[1,kd]\backslash[(k_0-1)d+1,(k_0+1)d]}MP^*_{[1,kd]\backslash\{\gamma\cup\{j\}\}}=\begin{pmatrix}B_1&0\\ B_2&B_3
\end{pmatrix}
$$
where the first column of $B_3$ is zero. Hence
$$
\det\left({P_{[1,kd]\backslash[(k_0-1)d+1,(k_0+1)d]}MP^*_{[1,kd]\backslash\{\gamma\cup\{j\}\}}}\right)=\det{B_1}\det{B_3}=0.
$$
\end{pf}
By Proposition \ref{ppp2}, we have
$$
\#\left\{\{\gamma_j\}_{j=1}^{2d-1}\cap[(k_0-2)d+1,k_0d]\right\}=d,$$
$$
\#\left\{\{\gamma_j\}_{j=1}^{2d-1}\cap[k_0d+1,(k_0+2)d]\right\}=d-1.
$$
Hence
\begin{align*}
|M(i,j)|\leq &\sum\limits_{\gamma}\left|\det\left({P_{[(k_0-1)d+1,(k_0+1)d]\backslash\{i\}}MP^*_{\gamma}}\right)\det\left({P_{[1,kd]\backslash[(k_0-1)d+1,(k_0+1)d]}MP^*_{[1,kd]\backslash\{\gamma\cup\{j\}\}}}\right)\right|\\
=&\sum\limits_{\sigma}\sum\limits_{\tau}|\det\left({P_{[(k_0-1)d+1,(k_0+1)d]\backslash\{i\}}MP^*_{\sigma\cup\tau}}\right)\mu_{\sigma}(i,j)\mu^{\tau}(i,j)|.
\end{align*}
\end{pf}

\section{The lower bound of the averaged denominator of the Green's function}
Recall that
$$
P_n(V,\theta,E)=\det{R_{[0,n-1]}(L_{V,\alpha,\theta}-E)R^*_{[0,n-1]}}.
$$
Based on a combination of the ideas in \cite{ajs} and \cite{gjyz}, we prove that for $n$ sufficiently large, we have
\begin{equation}\label{mai1}
\frac{1}{n}\int_{\T}\ln|P_{n}(V,\theta,E)|d\theta\geq \sum\limits_{j=1}^d\gamma_j(E)+\ln|V_d|+o(1).
\end{equation}
The proof involves complexification, first used by Herman \cite{H}. Note that by Jensen's formula, for any $\e>0$, one has
\begin{equation}\label{mai3}
\frac{1}{n}\int_{\T}\ln|P_{n}(V,\theta+i\e,E)|d\theta\leq \frac{1}{n}\int_{\T}\ln|P_{n}(V,\theta,E)|d\theta+2\pi\e.
\end{equation}
This follows from the fact $P_n(V,\theta,E)$ is a real trigonometric polynomial of degree $n$ as a function of $\cos2\pi(\theta)$.

On the other hand, if $E$ is subcritical, by Theorem \ref{sub}, $\gamma_d(E)$ (the minimal Lyapunov exponent of $(\alpha,L_{E,V}^{2\cos})$) is larger than $0$. Thus there is a gap between the $d$-th and $(d+1)$-th Lyapunov exponent of $(\alpha,L_{E,V}^{2\cos})$. Following the proof of Theorem 1.4 in \cite{ajs} (see also footnote 17 in \cite{ajs}),
$$
\omega^d(E)=\lim\limits_{\e\rightarrow 0^+}\frac{1}{2\pi\e}(L^d(L_{E,V}^{2\cos}(\cdot+\e))-L^d(L_{E,V}^{2\cos}))
$$
is an integer. It is also proved in \cite{gjyz} (see Proposition 5.1) that $|\omega^d(E)|\leq1$. Thus $\omega^d(E)=1$ \footnote{Note that $\omega^d(E)>0$ if $E\in \Sigma_{V,\alpha}$.}. It follows that
\begin{equation}\label{mai2}
L^d(L_{E,V}^{2\cos}(\cdot+\e))=L^d(L_{E,V}^{2\cos})+2\pi\e=\sum\limits_{j=1}^d\gamma_j(E)+2\pi\e
\end{equation}
for $\e$ sufficiently small. \eqref{mai1} follows immediately from \eqref{mai3}, \eqref{mai2} and the following equality (for subcritical $E$ and sufficiently small $\e$ with $\e\neq 0$),
\begin{equation}\label{mai4}
\lim\limits_{n\rightarrow\infty} \frac{1}{n}\int_{\T}\ln|P_{n}(V,\theta+i\e,E)|d\theta=L^d(L_{E,V}^{2\cos}(\cdot+\e))+\ln|V_d|.
\end{equation}
We remark, \eqref{mai4} is hard to prove when $\e=0$, while the case $\e\neq 0$ is much more simple, the main idea here is to involve dominated splitting.

Note that since $L^d(L_{E,V}^{2\cos}(\cdot+i\e))$ is a piecewise affine function (see \cite{ajs}), thus for $\e$ sufficiently small, $(\alpha,L_{E,V}^{2\cos}(\cdot+i\e))$ is $d$-regular, hence $d$-dominated. By the definition of domination (we now fixed $E$ and $\e$), one has
\begin{Corollary}\label{m'}
There are sequences of $d\times d$ matrix valued functions $\{F_{\pm}(k,\theta+i\e,E)\}_{k\in\Z}\in C^\omega(\T, \rm{M}(d,\C))$ obeying
$$
C^*F_{\pm}(k-1,\theta+i\e,E)+CF_\pm(k+1,\theta+i\e,E)+B(\theta+i\e+kd\alpha)F_\pm(k,\theta+i\e,E)=EF_\pm(k,\theta+i\e,E),
$$
$$
\sum\limits_{k=0}^\infty\|F_+(k,\theta+i\e,E)\|^2<\infty, \ \ \sum\limits_{k=-\infty}^{0}\|F_-(k,\theta+i\e,E)\|^2<\infty.
$$
\end{Corollary}
By invariance, there are $\Lambda^\pm(\theta+i\e,E)$ such that
\begin{align*}
&L_{E,V}^{2\cos}(\theta+i\e)\begin{pmatrix}F^+(1,\theta+i\e,E)&F^-(1,\theta+i\e,E)\\ F^+(0,\theta+i\e,E) &F^-(0,\theta+i\e,E)\end{pmatrix}\\
=&\begin{pmatrix}F^+(1,\theta+\alpha+i\e,E)&F^-(1,\theta+\alpha+i\e,E)\\ F^+(0,\theta+\alpha+i\e,E) &F^-(0,\theta+\alpha+i\e,E)\end{pmatrix}\begin{pmatrix}\Lambda^+(\theta+i\e,E)&0\\ 0&\Lambda^-(\theta+i\e,E)\end{pmatrix}.
\end{align*}
We denote  $\Theta_n=\left\{\theta:\det{F^\pm (0,\theta+n\alpha+i\e,E)}=0\right\}\cup \left\{\theta:\det{F^\pm (1,\theta+n\alpha+i\e,E)}=0\right\}$ and
\begin{align*}
&\left(L_{E,V}^{2\cos}\right)_n(\theta+i\e)=\begin{pmatrix} A^n_{11}(\theta+i\e,E)& A^n_{12}(\theta+i\e,E)\\ A^n_{21}(\theta+i\e,E)&A^n_{22}(\theta+i\e,E)\end{pmatrix}\\
=&\begin{pmatrix}F^+(1,\theta+i\e+n\alpha,E)&F^-(1,\theta+i\e+n\alpha,E)\\ F^+(0,\theta+i\e+n\alpha,E) &F^-(0,\theta+i\e+n\alpha,E)\end{pmatrix}\begin{pmatrix}\Lambda_n^+(\theta+i\e,E)&0\\ 0&\Lambda_n^-(\theta+i\e,E)\end{pmatrix}\\
&\cdot \begin{pmatrix}F^+(1,\theta+i\e,E)&F^-(1,\theta+i\e,E)\\ F^+(0,\theta+i\e,E) &F^-(0,\theta+i\e,E)\end{pmatrix}^{-1}
\end{align*}
where
$$
\Lambda_n^\pm(\theta+i\e,E)=\Lambda_n^\pm(\theta+i\e+(n-1)\alpha,E)\cdots \Lambda_n^\pm(\theta+i\e+\alpha,E)\Lambda_n^\pm(\theta+i\e,E).
$$
When $\theta\in \Theta_n\cap \Theta_0$, we have \footnote{we omit $E$ in the following formula for simplicity.}
\small \begin{align*}
&A_{11}^n(\theta+i\e)=\left(F^+(0,\theta+i\e+n\alpha)-F^-(0,\theta+i\e+n\alpha)(F^-(1,\theta+i\e+n\alpha))^{-1}F^+(1,\theta+i\e+n\alpha)\right)^{-1}\\
&\cdot\left(-F^-(0,\theta+i\e+n\alpha)(F^-(1,\theta+i\e+n\alpha))^{-1}\Lambda^+_n(\theta+i\e)F^+(1,\theta+i\e)(\Lambda^-_n(\theta+i\e)F^+(0,\theta+i\e))^{-1}+I_d\right)\\
&\cdot\Lambda^-_n(\theta+i\e)F^+(0,\theta+i\e).
\end{align*}
Thus
\small \begin{align*}
&\ln|\det{A_{11}^n(\theta+i\e)}|=-\ln\left|\det{\begin{pmatrix}F^+(1,\theta+i\e+n\alpha)&F^-(1,\theta+i\e+n\alpha)\\ F^+(0,\theta+i\e+n\alpha) &F^-(0,\theta+i\e+n\alpha)\end{pmatrix}}\right|+\ln|\det{F^-(1,\theta+i\e+n\alpha)}|\\
&+\ln|\det{\left(-F^-(0,\theta+i\e+n\alpha)(F^-(1,\theta+i\e+n\alpha))^{-1}\Lambda^+_n(\theta+i\e)F^+(1,\theta+i\e)(F^+(0,\theta+i\e))^{-1}(\Lambda^-_n(\theta+i\e))^{-1}+I_d\right)}|\\
&+\ln|\det{\Lambda_n^-(\theta+i\e)}|+\ln|\det{F^+(0,\theta+i\e)}|.
\end{align*}
Note that
$$
\|\Lambda^+_n(\theta+i\e)\|,\ \ \|(\Lambda^-_n(\theta+i\e))^{-1}\|\leq C\lambda^{-n},
$$
uniformly for any $\theta\in\T$.

On the other hand, by analyticity,
$$
|\widetilde{\Theta}_n|:=\left|\{\theta+i\e:|\det{F^\pm(1,\theta+i\e+n\alpha)}|\leq (\lambda/2)^{-n}\ \ {\rm or} \ \ |\det{F^\pm(0,\theta+i\e)}|\leq (\lambda/2)^{-n}\}\right|\leq C\lambda^{-Cn}.
$$
and for $\theta\notin \widetilde{\Theta}_n$, one has
$$
\|(F^\pm(1,\theta+i\e+n\alpha))^{-1}\|, \ \ \|(F^\pm(0,\theta+i\e))^{-1}\| \leq C(d)(\lambda/2)^n.
$$
Hence
$$
\lim\limits_{n\rightarrow \infty}\int_\T\frac{1}{n}\ln|\det{A_{11}^n(\theta+i\e,E)}|d\theta=\lim\limits_{n\rightarrow \infty}\int_\T\frac{1}{n}\ln|\det{\Lambda_n^-(\theta+i\e,E))}|d\theta=L^d(L_{E,V}^{2\cos}(\cdot+\e)).
$$
Note that
$$
V_d^{-n}\det{P_n(V,\theta+i\e,E)}=\det{A_{11}^n(\theta+i\e,E)},
$$
we thus finish the proof.
\section{almost localization}
Once we have the upper bound of the numerator and the lower bound of averaged denominator of the Green's function, we can establish the almost localization of the finite-range operator \eqref{fo}, following the idea of \cite{A1,AvilaJito,AvilaJito1,lyzz}, in a standard way.
\subsection{Statement of the main theorem}
\begin{definition}[Resonances]
Fix $\epsilon_0>0$ and $\theta\in \T$. An integer $k\in \Z$ is called an \textit{$\epsilon_0-$resonance} of $\theta$ if $\|2 \theta - k\alpha\|_{\T}\leq e^{-\epsilon_0 |k|}$ and $\|2 \theta - k\alpha\|_{\T}=\min_{|l|\leq |k|} \|2 \theta - l\alpha\|_{\T}$. We denote by $\{n_l\}_l$  the set  of $\epsilon_0-$resonances of $\theta$, ordered in such a way that $|n_1|\leq |n_2|\leq \cdots$.
We say that $\theta$ is  $\epsilon_0-$\textit{resonant} if the set $\{n_l\}_l$ is infinite.
\end{definition}

\begin{definition}[Almost localization]
The family $\{L_{V,\alpha,\theta}\}_{\theta\in\T}$ is said to be \textit{almost localized} for $E\in\R$ if there exist constants $C_0, C_1,\epsilon_0, \epsilon_1>0$ such that for all $\theta\in\T$, any generalized solution $u=(u_k)_{k \in \Z}$ to the eigenvalue problem $L_{V,\alpha,\theta} u = E u$ with $u_0=1$ and $|u_k| \leq 1+|k|$ satisfies
\begin{equation}\label{almost_localization}
|u_k| \leq C_1 e^{-\epsilon_1 |k|},\quad \forall \ C_0 |n_j| \leq |k| \leq C_0^{-1} |n_{j+1}|,
\end{equation}
where $\{n_l\}_l$  is the set  of  $\epsilon_0-$resonances of $\theta$.
\end{definition}
Our main purpose in this section is to prove the following Theorem.
\begin{Theorem}\label{thm_almost_almost-1}
Let $\alpha \in \R\backslash \Q$ satisfy $\beta(\alpha)=0$. If $E \in \Sigma_{V,\alpha}$ \footnote{$\Sigma_{V,\alpha}$ is the spectrum of $H_{V,\alpha,\theta}$.} with $\gamma_d(E)>0$, then $\{L_{V,\alpha,\theta}\}_{\theta\in\T}$ is almost localized for $E$. Moreover, for any $\delta\in (0,\gamma_d(E))$, any $C_0>1$, there exists $\epsilon_0=\epsilon_0(V,C_0,\delta)>0$ such that the following holds. Let  $L_{V,\alpha,\theta} u =E u$, with $|u_j|\leq 1$ for all $j \in \Z$.
 \begin{enumerate}
\item  If $\theta$ is not $\epsilon_0-$resonant, then
$|u_{j}|\leq e^{-(\gamma_d(E) -\delta)|j|}$ for $|j|$ large enough.
\item Else, let $\{n_l\}_l$ be the set of $\epsilon_0-$resonances of $\theta$. Given any $\eta>0$,
\begin{equation*}
|u_{j}|\leq e^{-(\gamma_d(E)-\delta)|j|},\quad \forall \ 2C_0 |n_l|+\eta|n_{l+1}|< |j| < (2C_0)^{-1} |n_{l+1}|,
\end{equation*}
provided that $|j|$ is large enough.
\end{enumerate}
\end{Theorem}
\subsection{Proof of Theorem \ref{thm_almost_almost-1}}
Given $\theta \in \T$ and $\epsilon_0>0$, we denote by $\{n_l\}_l$  the  set of $\epsilon_0-$resonances of $\theta$, i.e.,
$$
\|2 \theta - n_l \alpha\|_\T \leq e^{-\epsilon_0|n_l|},\quad \text{and}\quad \|2 \theta - n_l \alpha\|_\T=\min_{|m|\leq |n_l|} \|2 \theta - m \alpha\|_\T.
$$
Fix $\theta \in \R$, and let $u=(u_j)_{j\in\Z}$ be a  solution to $L_{V,\alpha,\theta}  u= E  u$, with $ u_0=1$ and $| u_{j}|\leq 1$ for all $j \in \Z$.
%The explicit estimates on the decay rate of $u$
%come from  the proof of almost localization in \cite{AvilaJito}.
Given an interval $I=[x_1,x_2] \subset \Z$ of length $N\geq 0$, recall that for any $x \in I$, we have
\begin{align}\label{exrpe green}
u_x=-\sum\limits_{y=x_1}^{x_1+d-1}\sum\limits_{k=y+1-x_1}^dG_I(x,y)V_{k}u_{y-k}-\sum\limits_{y=x_2-d+1}^{x_2}\sum\limits_{k=-d}^{y-x_2-1}G_I(x,y)V_{k}u_{y-k}.
\end{align}
Then by Cramer's rule, we have
\begin{equation*}\label{cramerrule}
G_I(x,y)=\frac{\mu_{x,y}}{P_{N}(V,\theta+x_1\alpha,E)}.
\end{equation*}
Given $\xi > 0$ and $m\in \N$, we say that $y \in \Z$ is $(\xi,m)-$\textit{regular} if there exists an interval $J=[x_1,x_2]\subset \Z$ of length $m$ such that $y \in J$ and
$$
|G_{J}(y,x_1+j)| < e^{-\xi |y-x_i|},\quad |y-x_1| \geq \frac{1}{7}m,\quad j=0,\cdots,d-1,
$$
$$
|G_{J}(y,x_2-j)| < e^{-\xi |y-x_i|},\quad |y-x_2| \geq \frac{1}{7}m,\quad j=0,\cdots,d-1.
$$

Let $(q_i)_{i \geq 1}$ be the sequence of denominators of best approximants of $\alpha$. We associate with any integer $C_0 |n_l| < |j| < C_0^{-1} |n_{l+1}|$ scales $\ell \geq 0$ and $s\geq 1$ so that
$$
2 s q_\ell \leq\zeta j <   \min(2(s+1)q_\ell,2q_{\ell+1}),
$$
where $\zeta:=\frac{1}{32}$ if $2|n_l|< j < 2^{-1} |n_{l+1}|$, and $\zeta:=\frac{C_0-1}{16C_0}$ otherwise.
We set
\begin{itemize}
\item $I_1:=[-2s q_\ell+1,0]$ and $I_2:=[j-2 s q_\ell +1,j+2 sq_\ell]$ if $j<|n_{l+1}|/3$, $n_{l}\geq 0$.
\item $I_1:=[1,2s q_\ell]$ and $I_2:=[j-2 s q_\ell +1,j+2 sq_\ell]$ if $j<|n_{l+1}|/3$ and $n_{l}<0$.
\item $I_1:=[-2s q_\ell+1,2s q_\ell]$ and $I_2:=[j-2 s q_\ell +1,j]$ if $|n_{l+1}|/3\leq j< |n_{l+1}|/2$.
\item $I_1:=[-2s q_\ell+1,2s q_\ell]$ and $I_2:=[j+1,j+2 sq_\ell]$ if $j\geq |n_{l+1}|/2$.
\end{itemize} In particular, the total number of elements in $I_1 \cup I_2$ is $6 sq_\ell$.
Fix $\delta>0$ arbitrary. If $\epsilon_0>0$ is chosen sufficiently small, %and since $\theta\equiv \la n \ra\ \mathrm{mod}\ \Z/2$ is non-resonant,
then in view of $\beta(\alpha)=0$, Lemma 5.8 in \cite{AvilaJito}  implies that there exists an integer $j_0=j_0(C_0,\alpha,\delta)>0$ such that for $j> j_0$, the set $\{\theta_m:=\theta+m \alpha\}_{m \in I_1 \cup I_2}$ is $\delta-$\textit{uniform}, i.e.,
$$
\max\limits_{z \in [-1,1]}\max\limits_{m \in I_1 \cup I_2} \prod\limits_{m\neq p\in I_1 \cup I_2} \frac{|z-\cos (2 \pi \theta_p)|}{|\cos (2 \pi \theta_m)-\cos(2 \pi \theta_p)|}< e^{(6sq_\ell-1)\delta}.
$$
Following the proof of Lemma 5.4 in \cite{AvilaJito} with the lower bound obtained in Section 5, we conclude that
%By Lemma 5.3 in \cite{AvilaJito}, for any $0<\delta<1/2$, $\epsilon>0$, there exist $N_0=N_0(\lambda,\alpha,\delta) \in \N$ and $\sigma=(\lambda,\alpha,\delta)>0$ such that if $y \in \Z$ is $(-\ln |\lambda|-\epsilon,N)-$singular for some $N>N_0$, then for any $x \in \Z$ such that $y-(1-\delta)N \leq x \leq y - \delta N$, we have $|P_N(\theta+x \alpha)|\leq e^{(N+1)(-\ln |\lambda|-\sigma)}$.
for any $\delta>0$, there exists $j_1=j_1(C_0,\alpha,V,d,\delta)>0$ such that any $j> j_1$ is $(\gamma_d(E)-\delta,6 s q_\ell-1)-$regular.\\
\textbf{Proof of Theorem \ref{thm_almost_almost-1} (2)}:
We consider the case that $\theta$ is $\epsilon_0-$resonant.
We will show that the sequence $( u_j)_j$ decays exponentially in some suitable interval between two consecutive resonances, with a rate close to the $d$-th Lyapunov exponent $\gamma_d(E)$.
By the condition $\beta(\alpha)=0$, we know that $|n_l|=o(|n_{l+1}|)$. Let us fix some small $\delta>0$. Given $l>0$ sufficiently large, take $\ell>0$ such that $2q_\ell \leq \zeta(2C_0 |n_l|+1)<2q_{\ell+1}$, and let $2C_0 |n_l|+\eta|n_{l+1}|\leq |j| \leq (2C_0)^{-1} |n_{l+1}|$. We set $b_l:=2 C_0 |n_l|+1$.
Then
for any $y \in [b_l,2j]$, there exists an interval $I(y)=[x_1,x_2]\subset \Z$ %\subset [-4j,4j]$
%of length $6sq_\ell-1$
with $y \in I(y)$ and
$$
\mathrm{dist}(y,x_i)\geq \frac{1}{7}|I(y)|\geq \frac{6 q_{\ell}-1}{7}\geq \frac{q_{\ell}}{2},
$$
such that
$$
|G_{I(y)}(y,x_1+j)| < e^{-(\gamma_d(E)-\delta) |y-x_i|}\leq e^{-(\gamma_d(E)-\delta) \frac{q_\ell}{2}},\quad j=0,\cdots,d-1.
$$
$$
|G_{I(y)}(y,x_2-j)| < e^{-(\gamma_d(E)-\delta) |y-x_i|}\leq e^{-(\gamma_d(E)-\delta) \frac{q_\ell}{2}},\quad j=0,\cdots,d-1.
$$
We denote by
$$
\partial I(y)=\{x_1-d+1,\cdots,x_1,x_2,\cdots,x_2+d-1\}.
$$
For $z \in \partial I(y)$,  we denote by
$$
z'=\begin{cases}\{j\}_{j=z+1-x_1}^d&z=x_1-d+1,\cdots,x_1\\
\{j\}_{j=-d}^{z-x_2-1}&z=x_2,\cdots,x_2+d-1
\end{cases}.
$$
If $x_2+d-1<2 j$ or $x_1-d+1>b_l$, we can expand $ u_{x_2+j}$ or $u_{x_1-j}$ for $j=0,\cdots,d-1$ following \eqref{exrpe green}, with $I=I(x_2+j)$ or $I=I(x_1-j)$. We continue to expand each term until we arrive to $\widetilde z$ such that either $\widetilde z\leq b_l$, or  $\widetilde z> 2j$, or the number %$s+1$
of $G_I$ terms in the following product becomes $\lfloor \frac{2j+2d}{q_{\ell}}\rfloor$, whichever comes first:
$$
 u_j=\sum\limits_{r,\ z_{i+1} \in \partial I(z_i')} G_{I(j)}(j,z_1)V_{z_1-z_1'}G_{I(z_1')}(z_1',z_2)\dots V_{z_r-z_r'}G_{I(z_r')}(z_r',z_{r+1}) u_{z_{r+1}'}.
$$
In the first two cases, we estimate
\begin{align*}
&|G_{I(j)}(j,z_1)V_{z_1-z_1'}G_{I(z_1')}(z_1',z_2)\dots V_{z_r-z_r'}G_{I(z_r')}(z_r',z_{r+1}) u_{z_{r+1}'}|\\ \leq\ &e^{-(\gamma_d(E)-\delta)(|j-z_{1}|+\sum_{i=1}^r |z_i'-z_{i+1}|)}\leq C(d)^{\frac{2j}{q_\ell}}e^{-(\gamma_d(E)-\delta)(|j-z_{r+1}|-(r+1)d)} \\
\leq\ &\max\big(C(d)^{\frac{2j}{q_\ell}}e^{-(\gamma_d(E)-\delta)(j-b_l-\frac{2j}{q_{\ell}})},C(d)^{\frac{2j}{q_\ell}}e^{-(\gamma_d(E)-\delta)(2j-j-\frac{2j}{q_{\ell}})}\big) \leq e^{-(\gamma_d(E)-\delta)(j+o(j))},
\end{align*}
where we have used that $|b_l| = o(|j|)$, while in the third case, we have
\begin{align*}
&|G_{I(j)}(j,z_1)V_{z_1-z_1'}G_{I(z_1')}(z_1',z_2)\dots V_{z_r-z_r'}G_{I(z_r')}(z_r',z_{r+1}) u_{z_{r+1}'}\leq C(d)^{\frac{2j}{q_\ell}}e^{-(\gamma_d(E)-\delta)\frac{q_{\ell}}{2}\lceil \frac{2j+d}{q_{\ell}}\rceil}.
\end{align*}
Fix $\delta>0$ arbitrarily small.
By taking $|j|$ to be sufficiently large, resp. $\eta$ small enough in the previous expression, we conclude that $| u_j| \leq e^{-(\gamma_d(E)-\delta) |j|}$ for $|j|$ large enough with $2C_0 |n_l|+\eta|n_{l+1}| \leq |j| \leq (2C_0)^{-1} |n_{l+1}|$.\qed \\
\textbf{Proof of Theorem \ref{thm_almost_almost-1} (1)}:
We consider the other case, i.e., when $\theta$ is not $\epsilon_0-$resonant. Denote by $n$ its last $\epsilon_0-$resonance, set $b:=2 C_0 |n|+1$ and let $|j | \geq b$. Let us fix some small $\eta>0$.
Then
for any $y \in [b,2j]$, there exists an interval $I(y)=[x_1,x_2]\subset \Z$ %\subset [-4j,4j]$
%of length $6sq_\ell-1$
with $y \in I(y)$ and
$$
\mathrm{dist}(y,x_i)\geq \frac{1}{7}|I(y)|\geq \frac{6 q_{\ell}-1}{7}\geq \frac{q_{\ell}}{2},
$$
such that
such that
$$
|G_{I(y)}(y,x_1+j)| < e^{-(\gamma_d(E)-\delta) |y-x_i|}\leq e^{-(\gamma_d(E)-\delta) \frac{q_\ell}{2}},\quad j=0,\cdots,d-1,
$$
$$
|G_{I(y)}(y,x_2-j)| < e^{-(\gamma_d(E)-\delta) |y-x_i|}\leq e^{-(\gamma_d(E)-\delta) \frac{q_\ell}{2}},\quad j=0,\cdots,d-1.
$$
As previously, we can expand $ u_{x_2+j}$ or $u_{x_1-j}$ following \eqref{exrpe green}, with $I=I(x_2+j)$ or $I=I(x_1-j)$ for $j=0,\cdots,d-1$. We continue to expand each term until we arrive to $\widetilde z$ such that either $\widetilde z\leq b$, or  $\widetilde z> 2j$, or the number %$s+1$
of $G_I$ terms in the following product becomes $\lfloor \frac{2j+2d}{q_{\ell}}\rfloor$, whichever comes first:
$$
 u_j=\sum\limits_{r,\ z_{i+1} \in \partial I(z_i')} G_{I(j)}(j,z_1)V_{z_1-z_1'}G_{I(z_1')}(z_1',z_2)\dots V_{z_r-z_r'}G_{I(z_r')}(z_r',z_{r+1}) u_{z_{r+1}'}.
$$
In the first two cases, we estimate
\begin{align*}
&|G_{I(j)}(j,z_1)V_{z_1-z_1'}G_{I(z_1')}(z_1',z_2)\dots V_{z_r-z_r'}G_{I(z_r')}(z_r',z_{r+1}) u_{z_{r+1}'}|\\ \leq\ &e^{-(\gamma(E)-\delta)(|j-z_{1}|+\sum_{i=1}^r |z_i'-z_{i+1}|)}\leq C(d)^{\frac{2j}{q_\ell}}e^{-(\gamma_d(E)-\delta)(|j-z_{r+1}|-(r+1)d)} \\
\leq\ &\max\big(C(d)^{\frac{2j}{q_\ell}}e^{-(\gamma_d(E)-\delta)(j-b_l-\frac{2j}{q_{\ell}})},C(d)^{\frac{2j}{q_\ell}}e^{-(\gamma_d(E)-\delta)(2j-j-\frac{2j}{q_{\ell}})}\big)\leq e^{-(\gamma_d(E)-\delta)(j+o(j))},
\end{align*}
where we have used that $|b_l| = o(|j|)$, while in the third case, we have
\begin{align*}
&|G_{I(j)}(j,z_1)V_{z_1-z_1'}G_{I(z_1')}(z_1',z_2)\dots V_{z_r-z_r'}G_{I(z_r')}(z_r',z_{r+1}) u_{z_{r+1}'}|\leq C(d)^{\frac{2j}{q_\ell}}e^{-(\gamma_d(E)-\delta)\frac{q_{\ell}}{2}\lceil \frac{2j+d}{q_{\ell}}\rceil}.
\end{align*}
Fix $\delta>0$ arbitrarily small, by taking $|j|$ to be sufficiently large, we conclude that  $| u_j| \leq e^{-(\gamma_d(E)-\delta) |j|}$ for $|j|$ large enough.
\qed

\section{Almost reducibility}
Almost localization of the finite-range operator \eqref{fo} immediately implies almost reducibility of the original Schr\"odinger cocycle also in a standard way.
\subsection{Statement of the main result}
\begin{Theorem}\label{alal}
Let $\alpha \in \R\backslash\Q$ satisfy $\beta(\alpha)=0$ and $E\in \Sigma_{V,\alpha}$. If $(\alpha, A_E)$ is subcritical with subcritical radius $h(E)>0$, then for any $r\in (0,h)$, the following holds  on $\{|\Im z|<r\}$:
\begin{enumerate}
\item

either $(\alpha,A_E)$ is almost reducible to $(\alpha,R_{\theta})$
 for some
$\theta=\theta(E)\in \R$:  for any $\varepsilon>0$, there exists $U\in C^\omega_{r}(\T,\mathrm{PSL}(2,\R))$ such that
\begin{equation*}
|U(\cdot+\alpha)^{-1} A_E(\cdot) U(\cdot)-R_{ \theta}|_{r} < \varepsilon;
\end{equation*}

\item or $(\alpha,A_E)$ is reducible:
\begin{enumerate}
\item if $2 \rho(E)-j\alpha \notin \Z$ for any $j \in \Z$, then $(\alpha,A_E)$ is reducible to $(\alpha,R_\theta)$ for some
$\theta=\theta(E)\in \R$;
\item if $2 \rho(E)-k\alpha \in \Z$ for some $k \in \Z$, then $(\alpha,A_E)$ is reducible to
$$
\begin{pmatrix}
1 & c \\
0 & 1
\end{pmatrix}.
$$
\end{enumerate}
\end{enumerate}
\end{Theorem}
\subsection{Complex almost reducibility}

Fix $0<\tilde{r}<h$ and $\delta>0$, set $r:=\frac{1}{2}(h+\tilde r)$, $\eta:=\frac{h-r}{1000C_0}$. For any $E\in\Sigma_{V,\alpha}$, there exist $\theta=\theta(E)$ and $\{u_{j}\}_{j\in\mathbb{Z}}$ satisfying $L_{V,\theta,\alpha}u=Eu$ with $u_{0}=1$ and $|u_{j}|\leq 1$ for every $j\in\mathbb{Z}$. We consider the case that $\theta$ is $\epsilon_{0}$-resonant. Assume $\{n_{l}\}$ are resonances of $\theta$. For $l$ sufficiently large, we have
\begin{equation}\label{ales}
|u_{j}|\leq e^{-(\gamma_d -\delta)|j|},\quad \forall \ 2C_0 |n_l|+\eta|n_{l+1}| < |j| < (2C_0)^{-1} |n_{l+1}|,
\end{equation}

Denote $n:=|n_{l}|+1$, $N:=|n_{l+1}|$. Let $I:=[-[\frac{N}{C_0}]+1,[\frac{N}{C_0}]-1]=[x_{1},x_{2}]$. We define $u^{I}:z\rightarrow\Sigma_{j\in I}u_{j}e^{2\pi ijz}$ and $U^I:z\rightarrow\begin{pmatrix} e^{2\pi i\theta}u^{I}(z) \\ u^{I}(z-\alpha)\end{pmatrix}$. Then
\begin{align}\label{eq}
A_E(z)U^I(z)=e^{2\pi i\theta}U^I(z+\alpha)+\begin{pmatrix} g(z) \\ 0 \end{pmatrix},
\end{align}
for some function $g\in C^\omega(\T,\C)$ whose Fourier coefficients $(\widehat{g}_j)_{j \in \mathbb{Z}}$ satisfy
$$
\widehat g_j=\chi_{I}(j) \left(E-2 \cos 2\pi(\theta+j\a)\right) u_j - \sum_{k=-d}^{d} \chi_{I}(j-k)  u_{j-k} V_k.
$$
Since $L_{V,\theta,\alpha} u=Eu$, we also have
\begin{equation}\label{equation coeffc gg}
\widehat g_j=-\chi_{\Z \backslash I}(j) \left(E-2 \cos 2\pi(\theta+j\a)\right)  u_j + \sum_{k=-d}^d \chi_{\Z \backslash I}(j-k) u_{j-k} V_k.
\end{equation}
If $E \in \Sigma_{V,\alpha}$, then $|E|\leq 2+ |V|_0 \leq C(V)$. By \eqref{ales} and \eqref{equation coeffc gg}, we therefore obtain
%for any $0\leq \delta< h_1$, there exist $\widetilde C_i=\widetilde C_i(C_1,\epsilon_1,h_0,\delta)>0$, $i=3,4,5$, such that
%some constants $0<h'=h'(\epsilon_1,V)<h/2$ and $c'=c'(\epsilon_1,V)>0$,
\begin{align*}
|g|_{r} \leq& \ \sum_{j\in \Z} \chi_{\Z \backslash I}(j)C(d,V) e^{-(\gamma_d-\delta-r)|j|}+\sum_{j\in \Z}\sum_{k=-d}^d  \chi_{\Z \backslash I}(j-k)C(d,V) e^{-(\gamma_d-\delta)|j-k|} e^{- \gamma_d|k|}e^{r|j|}\\
\leq&\ C(d,V) e^{-(\gamma_d-\delta-r)C_0^{-1}|N|}.
\end{align*}

\begin{Lemma}[\cite{AvilaJito}]\label{gui}
Let $l\geq 1$ and $1\leq p\leq[q_{l+1}/q_{l}]$. If $P$ has essential degree at most $pq_{l}-1$ and $x_{0}\in\mathbb{T}$, then for some absolute constant $K_{0}>0$,
$$
\|P\|_{C^{0}}\leq K_{0}q_{l+1}^{K_{0}p}\sup_{0\leq m\leq pq_{l}-1}|P(x_{0}+m\alpha)|.
$$
\end{Lemma}
\begin{Lemma}\label{pes2}
We have
$$
e^{-2\eta |n_{l+1}|}\leq \inf_{|\Im z|\leq r}|U^I(z)|\leq \sup_{|\Im z|\leq r}|U^I(z)|\leq e^{2\eta |n_{l+1}|},
$$
for $l$ sufficiently large depending on $d$.
\end{Lemma}
\begin{pf}
We proof the first inequality by contradiction. Suppose for some $z_{0}\in\{|\Im z|\leq r\}$ we have $\|U(z_{0})\|<e^{-2\eta m}$.
Note that for any $l\in\mathbb{N}$,
\begin{align*}
e^{2\pi il\theta}U^{I}(z_{0}+n\alpha)=(A_E)_{n}(z_{0})U^{I}(z_{0})-\sum_{s=1}^{n}e^{2\pi i(s-1)\theta}(A_E)_{n-s}(z_{0}+s\alpha)\begin{pmatrix}g(z_{0}+(s-1)\alpha)\\0\end{pmatrix}.
\end{align*}
Since $L_\e(E)=0$ for $|\e|\leq\gamma_d$, by similar arguments as Corollary 2.1 in \cite{lyzz}, we have
$$
\sup_{|\Im z|\leq r}\|(A_E)_n(z)\|\leq Ce^{\eta |n|}, \ \ \forall n\in\N.
$$
Thus
$$
\|U^{I}(z_{0}+n_{l+1}\alpha)\|\leq C|n_{l+1}|e^{\eta |n_{l+1}|}
e^{-(\gamma_d-r-\delta)C_0^{-1}|n_{l+1}|}\leq e^{-\eta|n_{l+1}|},
$$
for $l$ sufficiently large. By Lemma \ref{gui} and $\beta(\alpha)=0$, we have
$$
\|u^{I}(x+i\Im z_{0})\|_{C^{0}}\leq K_{0}e^{o(n_{l+1})}e^{-\eta|n_{l+1}|}<1.
$$
This contradicts to $u_0=1$.

On the other hand, by \eqref{ales}, we have
$$
\sup_{|\Im z|\leq r}|U^I(z)|\leq e^{2\eta |n_{l+1}|}.
$$
\end{pf}

\begin{Lemma}[\cite{A1}]\label{58}
Let $V:\mathbb{T}\rightarrow\mathbb{C}^{2}$ be analytic in $|\Im z|<\eta$. Assume that $\delta_{1}<\|V(z)\|<\delta_{2}^{-1}$ holds on $|\Im z|<\eta$. Then there exists $M:\mathbb{T}\rightarrow\mathrm{SL}(2,\mathbb{C})$ analytic on $|\Im z|<\eta$ with first column $V$ and $\|M\|_{\eta}\leq C\delta_{1}^{-2}\delta_{2}^{-1}(1-\ln(\delta_{1}\delta_{2}))$.
\end{Lemma}
\begin{Lemma}\label{lemma first conjugaison}
There exist constants $C=C(\alpha,d,V,r)>0$, such that the following holds. There exists $B\in C^{\omega}(\T, \mathrm{SL}(2,\C))$ with
$\|B\|_{r} \leq Ce^{4\delta |n_{l+1}|}$ such that
\begin{equation}\label{conj par complexe bis}
B(\cdot+\a)^{-1} A_E(\cdot) B(\cdot)=\begin{pmatrix}
e^{2 \pi \mathrm{i} \theta} & 0\\[1mm]
0 & e^{-2 \pi \mathrm{i} \theta}
\end{pmatrix}+\begin{pmatrix}
\beta_1(\cdot) & \beta_2(\cdot)\\[1mm]
\beta_3(\cdot) & \beta_4(\cdot)
\end{pmatrix},
\end{equation}
with  $|\beta_j|_{r} \leq C e^{-\delta|n_{l+1}|} $ for $j=1,2,3,4$.
\end{Lemma}
\begin{pf}
By Lemma \ref{gui} and Lemma \ref{58}, we can define $U_2 \in C^{\omega}_{r}(\T,\mathrm{SL}(2,\C))$ with first column $\mathcal{U}^J$ such that
\begin{equation*}
U_2(\cdot+\a)^{-1} A_E(\cdot) U_2(\cdot)=\begin{pmatrix}
e^{2 \pi \mathrm{i} \theta} & 0\\[1mm]
0 & e^{-2 \pi \mathrm{i} \theta}
\end{pmatrix}+\begin{pmatrix}
\widetilde \beta_1(\cdot) &b(\cdot)\\[1mm]
\widetilde\beta_3(\cdot) & \widetilde\beta_4(\cdot)
\end{pmatrix}
\end{equation*}
with $|U_2|_{r} \leq C(\alpha,d,V,r)  e^{3\delta  n_{l+1}}$, $|b|_{r} \leq C(\alpha,d,V,r) e^{2\delta n_{l+1}}$, and for $j=1,3,4$,
$$
|\widetilde\beta_j|_{r} \leq C(\alpha,d,V,r) e^{-(\gamma_d-\delta-r) C_0^{-2}n_{l+1}}.
$$

Let us write $b(z)=\sum_j\widehat b_j e^{2 \pi \mathrm{i} j z}$, and let  $\tau$ satisfy $$b (z)-e^{-2 \pi \mathrm{i} \theta}\tau(z+\alpha)+e^{2 \pi \mathrm{i} \theta} \tau(z)= \sum_{|j| \geq n_{l+1}} \widehat b_j e^{2 \pi \mathrm{i} j z}.$$ We have $\tau(z)=\sum_{|j| < n_{l+1}} \widehat \tau_j e^{2 \pi \mathrm{i} j z}$, where
$
\widehat\tau_j:=\frac{-\widehat b_j e^{-2 \pi \mathrm{i} \theta}}{1-e^{-2 \pi \mathrm{i} (2\theta-j\alpha)}}
%=\frac{- \widehat \varphi_j e^{- \pi \mathrm{i} n\alpha}}{1-e^{-2 \pi \mathrm{i} (n-j)\alpha}}
$.
By the assumption   $\beta(\alpha)=0$, and the definition of resonances, for $j \neq n_{l+1}$, we have
$$\|2\theta- j \alpha\|_{\T}  \geq \|(j-n_{l+1})\alpha\|_{\T}- \|2\theta- n_{l+1}\alpha\|_{\T} \geq  e^{-o(|j-n_{l+1}|)}-e^{-\epsilon_0n_{l+1} }\geq \frac{1}{2}e^{-o(|j-n_{l+1}|)}.$$
%\begin{align*}
%\|2\theta- j \alpha\|_{\T}  \geq\ & \|(j-n_l)\alpha\|_{\T}- \|2\theta- n_l\alpha\|_{\T} \\\geq\ &  e^{-o(|j-n_l|)}-e^{-\epsilon_0n_l }\geq \frac{1}{2}e^{-o(|j-n_l|)}.
%\end{align*}
Therefore,  we deduce that $|\tau|_{r} \leq C(\alpha,d,V,r)e^{2\delta n_{l+1}}$.

Let $B:=U_2
\begin{pmatrix}
1 & \tau\\
0 & 1
\end{pmatrix}$ conjugate the initial cocycle to the following:
$$
B(\cdot+\a)^{-1} A_E(\cdot) B(\cdot)=\begin{pmatrix}
e^{2 \pi \mathrm{i} \theta} & 0  \\[1mm]
0 & e^{-2 \pi \mathrm{i} \theta}
\end{pmatrix}+
\begin{pmatrix}
\beta_1(\cdot) & \beta_2(\cdot)+\varsigma(\cdot)\\[1mm]
\beta_3(\cdot) & \beta_4(\cdot)
\end{pmatrix}
$$
with $\varsigma:z\mapsto \sum_{|j| \geq n_{l+1}} \widehat b_{j} e^{2 \pi \mathrm{i} j z}$,
 $\beta_1(\cdot):=\widetilde{\beta}_1(\cdot)-\widetilde{\beta}_3(\cdot)\tau(\cdot+\alpha)$,  $\beta_3(\cdot):=\widetilde{\beta}_3(\cdot)$,  $\beta_4(\cdot):=\widetilde{\beta}_4(\cdot)+\widetilde{\beta}_3(\cdot)\tau(\cdot)$ and
$$
\beta_2(\cdot):=\widetilde{\beta}_1(\cdot) \tau(\cdot)- \widetilde{\beta}_4(\cdot)\tau(\cdot+\alpha)+\widetilde{\beta}_3(\cdot)\tau(\cdot)\tau(\cdot+\alpha).
$$
By the estimates on $\widetilde{\beta}_1,\widetilde{\beta}_3,\widetilde{\beta}_4$, we have $|\beta_j|_{r}\leq C(\alpha,d,V,r) e^{-2\delta n_{l+1}}$, for all $j=1,2,3,4$.
% Therefore, for all $r\in(0,\frac{\epsilon_1}{20 \pi})$, there exists $c_{10}=c_{10}(\alpha,h_0,r,\varsigma)>0$ such that $|\beta_j|_{4\varsigma}\leq c_{10} e^{-2 \pi c n_l}$.
On the other hand,
$$\big|\varsigma\big|_{\tilde r}\leq \sum_{|j| \geq n_{l+1}} |b|_{h}e^{-(r-\tilde r)|j|} \leq C(\alpha,d,V,r)e^{-2\delta n_l}.$$
\end{pf}
We immediately have the following Lemma.
\begin{Lemma}\label{upper}
For some $C>0$ depending on $\alpha,V$,
$$
\|(A_E)_n\|_{r}\leq C(1+|n|)^C.
$$
\end{Lemma}
\begin{pf}
Note that the above Lemmas work for any $4C_0(n_l+1)<m<C_0^{-1}|n_{l+1}|$ of the form $m=rq_k-1<q_{k+1}$. Let $n\geq C$, we can choose $\frac{\delta}{100}\ln n\leq m\leq 100\delta\ln n$ so that  $4C_0(n_l+1)<m<C_0^{-1}|n_{l+1}|$ and $m=rq_k-1<q_{k+1}$ for some $l$ and $k$. It follows from the following inequality
$$
\sup\limits_{0\leq n\leq e^{\frac{\delta}{100}m}}\|(A_E)_n\|_r\leq e^{8\delta m}.
$$
\end{pf}
\subsection{Real almost reducibility}
Let $\widetilde U(z)=e^{\pi in_lz}U(z)$ and $L^{-1}=\|2\theta-n_l\alpha\|$, we define $B(z)=(U(z),\overline{U(z)})$.
\begin{Lemma}\label{low}
We have
$$
\inf\limits_{x\in\T}|\det{B(x)}|\geq L^{-5C}.
$$
\end{Lemma}
\begin{pf}
Recall that for any $2\times 2$ complex matrix M with columns $V$ and $W$,
\begin{align}\label{dete}
|\det M|=\|V\|\min_{\lambda\in \mathbb{C}}\|W-\lambda V\|,
\end{align}
with the minimizing $\lambda$ satisfying $\|\lambda V\|\leq\|W\|$.

Now we prove the Lemma \ref{low} by contradiction. Suppose there exist some $x_{0}$ and $\lambda_{0}$ such that $\|e^{-\pi in_{l}x_{0}}\overline{\widetilde{U}(x_{0})}-\lambda_{0}e^{\pi in_{l}x_{0}}\widetilde{U}(x_{0})\|\leq L^{-4C}$. Note that
\begin{align*}
&\|e^{2\pi in\theta}U(x_{0}+n\alpha)-(A_E)_{n}(x_{0})\tilde{U}(x_{0})\|\leq \Sigma_{k=1}^{n}|(A_E)_{n-k}(x_{0}+k\alpha)\tilde{G}(x_{0}+(k-1)\alpha)|\leq C|n|^{C+1}e^{-2\delta n_{l+1}}.
\end{align*}
By Lemma \ref{upper}, For $0\leq n\leq L$, we have
\begin{align}\label{pes1}
\|e^{-2\pi in\theta}e^{-\pi in_{l}x_{0}}\overline{U(x_{0}+n\alpha)}-e^{2\pi in\theta}e^{\pi in_{l}x_{0}}\lambda_{0}U(x_{0}+n\alpha)\|\leq L^{-C}.
\end{align}
It follows
\begin{align*}
\|e^{-\pi in_{l}(x_{0}+n\alpha)}\overline{U(x_{0}+n\alpha)}-\lambda_{0}
e^{\pi in_{l}(x_{0}+n\alpha)}U(x_{0}+n\alpha)\| \leq  L^{-C}+2nL^{-1}\|U\|_{0}\leq L^{-\frac{1}{3}}, \quad 0\leq n\leq L^{\frac{1}{3}}.
\end{align*}
In the last inequality, we use the fact that $\|U\|_0\leq |n_l|^C$.

Let $m>(\ln L)^2$ be minimal of the form $m=q_k-1$. Set $J:=[-[m/2],m-[m/2]]$. Define $U^{J}:=\begin{pmatrix} e^{2\pi i\theta}u^{J} (x)\\ u^{J}(x-\alpha) \end{pmatrix}$. Since $\beta(\alpha)=0$, we have $m\leq C_0^{-1}n_{l+1}$. We have $\|U-U^{J}\|_{0}\leq L^{-1}$ by \eqref{ales}.

Using Lemma \ref{gui}, we have
\begin{align*}
&\|e^{-\pi in_{l}x}\overline{U^{J}(x)}-e^{\pi in_{l}x}\lambda_{0}U^{J}(x)\|_{0}\leq L^{-\frac{1}{3}}.
\end{align*}
Thus
\begin{align*}
&\|e^{-\pi in_{l}x}\overline{U(x)}-e^{\pi in_{l}x}\lambda_{0}U(x)\|_{0}\leq L^{-\frac{1}{4}}.
\end{align*}

Substituting $x_{1}=x_{0}+n\alpha$ in \eqref{pes1} and taking $n=[L/2]$, we get
$$
\|ie^{-\pi in_{l}x_{1}}\overline{U(x_{1})}+ie^{\pi in_{l}x_{1}}\lambda_{0}U(x_{1})\|\leq L^{-C}+2L^{-1}\|U\|_{0}.
$$
Thus we can get $U(x_{1})\leq L^{-\frac{1}{5}}$. This is a contradiction to the following fact,
$$
\inf\limits_{x\in\R/\Z}\|U\|_0\geq c|n_{l}|^{-C}.
$$
\end{pf}

\begin{Lemma}\label{lowdet}
Let $x_{0}\in \mathbb{T}$. For $l$ sufficiently large, we have
$$
\sup_{|\Im z|<\tilde{r}}|\det B(z)-\det B(x_{0})|\leq e^{-2\delta n_{l+1}}.
$$
\end{Lemma}
\begin{pf}
By construction of $B(z)$ and \eqref{eq} we have
$$
A_{E}(z)B(z)=\begin{pmatrix} e^{2\pi i(\theta-n_j\alpha)} & 0 \\ 0 & e^{-2\pi i(\theta-n_j\alpha)} \end{pmatrix}B(z+\alpha)+(G(z),\overline{G(\bar{z})}),
$$
where $\|G\|_{r}\leq C(\alpha,d,V,r)e^{-100\delta|n_{l+1}|}$. Thus we get
$$
|\det B(x_{0}+\alpha)-\det B(x_{0})|\leq 2\|\det B\|_{0}\|G\|_{0}\leq e^{-90\delta|n_{l+1}|}.
$$
This gives
$$
|\det B(x)-\det B(x_{0})|\leq e^{\delta|n_{l+1}|}100|n_{l+1}|e^{-90\delta|n_{l+1}|}\leq e^{-80\delta|n_{l+1}|},
$$
for $0\leq k\leq 100|n_{l+1}|$.

On the other have $\|U\|_r\leq e^{4\delta|n_{l+1}|}$. Thus
$$
|\det B(z)-\det B(x_{0})|\leq  e^{8\delta|n_{l+1}|}.
$$

Let $f(z)=\det B(z)-\det B(x_{0})$, by the Hadamard three-circle theorem,
$$
\ln\sup_{|\Im z|=\delta \frac{r+\tilde{r}}{2}}|f(z)|\leq(1-\delta)\ln\sup_{|\Im z|=0}|f(z)|+\delta\ln\sup_{|\Im z|=\frac{r+\tilde{r}}{2}}|f(z)|,\quad 0\leq\delta\leq\frac{r-\tilde{r}}{\tilde{r}+r}.
$$
We have $\|f\|_{\tilde{r}}\leq e^{-\delta|n_{l+1}|}$ .
\end{pf}
\subsection{Proof of Theorem \ref{alal}}
If $\theta$ is $\epsilon_0$-resonant, let $S:=\Re U$ and $T:-\Im U$. Let $W_{1}(x)$ be the matrix with columns $S$ and $T$. We have $B=W_{1}\begin{pmatrix} 1 & 1 \\ -i & i \end{pmatrix}$. Notice $\Re\det B(x)=0$ for $x\in \mathbb{T}$. Moreover, by Lemma \ref{lowdet}, we have
$$
\inf_{|\Im z|\leq r}|\Im \det B(z)|\geq L^{-5C}-e^{-2\delta n_{l+1}}\geq e^{-\delta n_{l+1}}.
$$
So $(\det W_{1})^{\frac{1}{2}}=(\frac{\det B}{2i})^{\frac{1}{2}}$ is analytic on $|\Im z|<r$. Let $W(x)=\frac{1}{(\det W_{1})^{\frac{1}{2}}}W_{1}$ we have
$$
\|W(\cdot+\alpha)^{-1}A_{E}W(\cdot)-R_{\theta}\|\leq e^{-\delta|n_{l+1}|}\leq \epsilon.
$$
with $W\in SL(2,\R)$ and
$$
\|W\|_r\leq e^{5\delta|n_{l+1}|}.
$$
If $\theta$ is $\epsilon_{0}$-non-resonant, by Theorem \ref{thm_almost_almost-1}, we have $|u_{j}|\leq e^{-2\pi \tilde{r}|j|}$ for $|j|>N$ where $N=N(\alpha,V,r)$. Define $u:z\rightarrow\Sigma_{j}u_{j}e^{2\pi ijz}$ and $U:z\rightarrow\begin{pmatrix} e^{2\pi i\theta}u(z) \\ u(z-\alpha) \end{pmatrix}$. Then
\begin{align}\label{pbl}
A_{E}(z)U(z)=e^{2\pi i\theta}U(z+\alpha),
\end{align}
Let $B$ be the matrix with columns $U(z)$ and $\overline{U(\bar{z})}$. So $B(z)$ is analytic in $|\Im z|\leq r$ and we have $\det B(x+\alpha)=\det B(x)$. We can get $\det B(z)$ is constant. If $\det B\neq 0$, let $W(x)=\frac{1}{(-2i\det B)}B\begin{pmatrix} 1 & i \\ 1 & -i \end{pmatrix}$, then we have
$$
W(x+\alpha)^{-1}A_E(x)W(x)=R_{\theta}.
$$

If $\det B=0$, In this case, we must have $2\theta=n\alpha(\mod\Z)$, by the minimality of $x \mapsto x+\alpha$, and the fact that $\widehat u_0=1$,  we see that $U$ does not vanish.
Define
${U}^{(1)}(z):=e^{ \pi {\rm i} n z} {U}(z)\in \C^2\backslash\{0\}$.
Since $2\theta -n\alpha \in \Z$, we get
\begin{equation}\label{eq widedtiel u}
A_{E}(\cdot) \, {U}^{(1)}(\cdot) = e^{\pi {\rm i}(2\theta -n\alpha)} \, {U}^{(1)}(\cdot+\alpha)= \pm \, {U}^{(1)}(\cdot+ \alpha). %\quad \forall x \in \T.
\end{equation}
Without loss of generality, we assume that $A_{E}(\cdot) \, {U}^{(1)} (\cdot)=   {U}^{(1)}(\cdot+ \alpha)$. There exist $\psi\colon \T \to \C$ with $|\psi|=1$ and $\mathcal{V}\colon \T \to \R^2\backslash\{0\}$ such that ${U^{(1)}}=\psi \, \mathcal{V}$ on $\T$.
By \eqref{eq widedtiel u}, we have
$$
A_{E} (x) \, \mathcal{V}(x)=\frac{\psi(x+\alpha)}{\psi(x)} \mathcal{V}(x+\alpha),\quad \forall \ x \in \T, $$ hence $\frac{\psi(x+\alpha)}{\psi(x)}\in \R$.
By the minimality of $x \mapsto x+\alpha$, we deduce that $\psi|_\T\equiv e^{2 \pi \mathrm{i} \theta_0}\in \C$ for some $\theta_0\in \R$. The map \begin{equation}\label{defv}\mathcal{V}\colon z \mapsto e^{ \pi \mathrm{i}(n z- 2 \theta_0)}\mathcal{U}(z)\end{equation} is analytic on $\{|\Im z|\leq r\}$ and satisfies
\begin{equation}\label{invarireim2}
A_{E}(z) \, \mathcal{V}(z) = \mathcal{V}(z+ \alpha),\quad   \forall  \ z\in \C/\Z \  {\rm with} \ |\Im z|\leq r.
\end{equation}
Applying Theorem \ref{58} to $\mathcal{V}$, we deduce that there exists % $\mathcal{W} \colon \mathbb{T}\to \C^2\backslash\{0\}$ such that the map
$U_1=(\mathcal{V}, \mathcal{W}) \in C_{r}^\omega(\mathbb{T}, \mathrm{SL}(2,\C))$ such that for all $r'\in(0,r]$,
\begin{equation}\label{rstim U1bisss}
 |U_1|_{r'} \leq C'(\alpha,V,d,r') e^{ (r'+\delta)|N|}.
\end{equation}
  Indeed, one can choose $U_1\in C_{r}^\omega(\mathbb{T}, \mathrm{PSL}(2,\R))$, since $\mathcal{V}|_{\mathbb{T}}$ takes values in $\R^2\backslash\{0\}$, then one only need to replace $U_1$ by $U_1=(\mathcal{V},\widetilde{\mathcal{W}})$, with $\widetilde{\mathcal{W}}\colon z \mapsto \frac{1}{2}(\mathcal{W}(z)+\overline{\mathcal{W}(\overline{z}}))$.

By \eqref{invarireim2}, there exists $\varphi^{(1)}\ \in C_{r}^\omega( \mathbb{T}, \R)$ such that \begin{equation}\label{eqatori}
U_1(\cdot+\a)^{-1} A_{E}(\cdot) U_1(\cdot)=
\begin{pmatrix}
1 & \varphi^{(1)}(\cdot)\\[1mm]
0 & 1
\end{pmatrix}.
\end{equation}
By \eqref{rstim U1bisss} and \eqref{eqatori}, we have:
$$|\varphi^{(1)}|_{r} \leq C(\alpha,d,V,r)e^{20 (r+\delta)|N|}.$$
Since $\beta(\a)=0$, we can solve the cohomological equation
\begin{equation}\label{cohomolllfognri}
\phi(z+\a) - \phi(z) =\varphi^{(1)}(z)-\int_\T \varphi^{(1)}(x) \, dx,
\end{equation}
with $\int_\T \phi(x)dx=0$. Moreover, $\phi\ \in C_{r}^\omega( \mathbb{T}, \R)$, and for any $r'\in(0,r]$, one has
\begin{equation}\label{esti_phi}
 |\phi|_{r'} \leq C(\alpha,d,V,r')e^{20(r'+\delta) |n|}.
\end{equation}
 Let  $U:=U_1
\begin{pmatrix}
1 & \phi\\
0 & 1
\end{pmatrix}$ and $\varphi:= \int_\T \varphi^{(1)}(x) \, dx$.
\eqref{eqatori} and \eqref{cohomolllfognri} implies that
\begin{equation}\label{equation conjugaison finale}
U(\cdot+\a)^{-1}S_{E}^{V}(\cdot)  U(\cdot)=
\begin{pmatrix}
1 & \varphi\\
0 & 1
\end{pmatrix} .
\end{equation}
Obviously, $U\in C_{r}^\omega(\mathbb{T}, \mathrm{PSL}(2,\R))$. By \eqref{rstim U1bisss} and \eqref{esti_phi}, we get the estimate
$$|U|_{r'} \leq C(\alpha,V,d,r') e^{40 r'|n|},\quad \forall \ r'\in(0, r].$$
Thus, we finish the proof.

\section{Acknowledgement}
The first version of this paper was finished in Oct, 2019 when I was a Visiting Assistant Professor at UCI and I would like to thank Professor Svetlana Jitomirskaya for many helpful discussions and kind support. We told Professor Artur Avila our different proof of the ARC in Dec 2019 and I would like to thank him for explaining to me his proof and the histories of the ARC. I would also like to thank Professor Jiangong You for many helpful discussions. I am also grateful to the organizers of QMath15,
Sep.2022, for giving me the opportunity to present some details of my
approach.

\end{document}